\documentclass[11pt]{elsarticle}
\usepackage{}
\usepackage{amssymb}
\usepackage{xcolor}
\usepackage{mathrsfs}
\usepackage{amscd}
\usepackage{appendix}
\usepackage{geometry}
\usepackage{setspace}
\usepackage{amsfonts}
\usepackage{amsmath}
\usepackage{lineno,hyperref}
\usepackage{tikz}
\usetikzlibrary{matrix,arrows}
\newtheorem{definition}{Definition}[section]
\newtheorem{theorem}[definition]{Theorem} 
\newtheorem{lemma}[definition]{Lemma}     

\numberwithin{equation}{section}

\biboptions{numbers,sort&compress} 
\journal{Journal of Differential Equations}
\begin{document}
\newgeometry{left=2.5cm,right=2.5cm,top=2.5cm,bottom=2.5cm}

\begin{frontmatter}

\title{Invariant curves of quasi-periodic reversible mappings and its application }

\author[]{Yan Zhuang}
\ead{yzhuangmath@163.com}

\author[]{Daxiong Piao}
\ead{dxpiao@ouc.edu.cn}

\author[]{Yanmin Niu \corref{cor1}}
\ead{niuyanmin@ouc.edu.cn}
\address{School of Mathematical Sciences,  Ocean University of China, Qingdao 266100, P. R. China}
\cortext[cor1]{Corresponding author}

\begin{abstract}
We consider the existence of invariant curves of real analytic reversible mappings which are quasi-periodic in the angle variables. By the normal form theorem, we prove that under some assumptions, the original mapping is changed into its linear part via an analytic convergent transformation, so that invariants curves are obtained. In the iterative process, by solving the modified homological equations, we ensure that the transformed mapping is still reversible. As an application, we investigate the invariant curves of a class of nonlinear resonant oscillators, with the Birkhoff constants of the corresponding Poincar$\acute{e}$ mapping all zeros or not.
\end{abstract}


\end{frontmatter}

\begin{spacing}{1.25}
\section{Introduction}
In this paper, we consider the mappings of the form
\begin{equation}\label{main-equa}
M:\theta_1=\theta+\gamma_0+f(\theta,r),\ \  r_1=r+g(\theta,r),
\end{equation}
where $f(\theta,r)$ and $g(\theta,r)$ are quasi-periodic in $\theta$ with frequency $\omega=(\omega_1,\omega_2,...,\omega_m)$, real analytic in a neighborhood of $r=0$ and $f(\theta,0)=g(\theta,0)=0$. Suppose mapping $M$ is reversible with respect to the involution $G:(x,y)\rightarrow(-x,y)$, that is, $GMG=M^{-1}$.

 Compared with the classical twist mappings
\begin{equation}\label{generaltwistcase}
\theta_1=\theta+\gamma_0+\sum_{k=0}^{+\infty}\gamma_kr^k+f(\theta,r),\ \  r_1=r+g(\theta,r),
\end{equation}
mapping $M$ in (\ref{main-equa}) can be seen as having no twist term with all Birkhoff constants $\gamma_1=\gamma_2=\cdots=0$. If some of the Birkhoff constants are not zero, $r^{-1}f,r^{-1}g$ are small and periodic in $\theta$, and mappings (\ref{generaltwistcase}) are area-preserving, Moser's twist theorem\cite{moserbook} tells us that in any neighborhood of the fixed point $r=0$, there are many analytic invariant closed curves surrounding this point, which implies the stability of the fixed point. When $f,g$ are quasi-periodic and the frequencies are sufficiently incommensurable with $2\gamma_0^{-1}\pi$, Zharnitsky in \cite{ZhV} obtained the same conclusion for exact symplectic map. In many applications, we may meet non-Hamiltonian systems, or systems have some symmetrical characteristics. Therefore, the theory for reversible system or for its corresponding Poincar\'{e} mapping, that is, reversible mapping has received widespread attention. The invariant curves theorem of reversible system was first obtained by Moser \cite{moser2} in 1965, and developed by himself \cite{moser3,moser4} and Sevryuk \cite{sevryuk} for both continuous and discrete systems. Based on the KAM technique, Liu \cite{Liu-1} proved the existence of invariant curves for quasi-periodic reversible mappings (\ref{generaltwistcase}), with $\gamma_1=1$ and the frequencies of $f,g$ and $2\gamma_0^{-1}\pi$ satisfying the Diophantine condition. Furthermore, within the framework of small twist theorem by Ortega\cite{ortega}, Liu extended the results of mapping with intersection property to the reversible mapping under periodic\cite{Liu-2} and quasi-periodic \cite{Liu-1} cases.

On the other hand, if all the Birkhoff constants are zero, there were several papers concerning this equation. In 2000, R\"{u}ssmann \cite{russman} studied an area-preserving mapping of the precise form
$$
x_1=x\cos\gamma_0-y\sin\gamma_0+\cdots,\ \  y_1=x\sin\gamma_0+y\cos\gamma_0+\cdots,
$$
where Birkhoff constants $\gamma_1=\gamma_2=\cdots=0$, and the dots denote the higher order terms in $x$ and $y$. He proved that if there is a formal change of variables transforming this mapping into the linearized normal form and $\gamma_0$ satisfies the Bruno condition, then there is a convergent change of variables taking the given mapping into its linear part. The stability of the fixed point $(0,0)$ therefore obtained. Recently, Hu and Liu et \cite{Liuhu} considered the existence of invariant curves of real analytic area-preserving mappings (\ref{generaltwistcase}) by the normal form method. They applied the results to the boundedness of all solutions for asymmetric oscillator.
As far as we know, there is few results about the existence of invariant curves for reversible mapping (\ref{main-equa}).

By R\"{u}ssmann's method and the assumption that a formal change takes (\ref{main-equa}) into the linearized form $(\theta,r)\rightarrow(\theta+\gamma_0,r)$, we aim to construct a convergent transformation $U$ which can be commutated with $G$, that is $GU=UG$, taking (\ref{main-equa}) into the linear part. In the follows, assume the constant $\gamma_0$  in mapping $M$ and frequency $\omega$ of functions $f,g$ satisfy the so-called Diophantine condition
\begin{equation}\label{dioph}
\displaystyle\left|\frac{\langle k,\omega\rangle\gamma_0}{2\pi}-j\right|\geq\frac{c_0}{|k|^\sigma},\quad k\in\mathbb{Z}^m\setminus\{0\},j\in\mathbb{Z},
\end{equation}
where $c_0$ is a small positive constant, $\sigma>0$. It not difficult to show that for $\sigma>m+1$, the Lebesgue measure for the set of $\gamma_0$ satisfying the above inequalities is positive for sufficient small $c_0$.

Now we can state our first main result.
\begin{theorem}\label{main-theo}
Assume that the mapping M in (\ref{main-equa}) is reversible with respect to the involution $G:(\theta,r)\rightarrow(-\theta,r)$, $\gamma_0$ satisfies (\ref{dioph}). If there exists a formal change of variables
\begin{equation}\label{form-chan}U:
\begin{cases}
\theta=\xi+u(\xi,\eta),\\
r=\eta+v(\xi,\eta),
\end{cases}
\end{equation}
where $u,v$ are quasi-periodic in $\xi$, such that (\ref{main-equa}) is transformed into the linearized form
\begin{equation}\label{norm-form}U^{-1}MU:
\begin{cases}
\xi_1=\xi+\gamma_0,\\
\eta_1=\eta,
\end{cases}
\end{equation}
then there exists a real analytic convergent change of variables  which transforms (\ref{main-equa}) into (\ref{norm-form}) and commutate with $G$.
\end{theorem}

Based on the classical KAM approach, we will transform (\ref{main-equa}) closer to the linear form (\ref{norm-form}) by a sequence of transformations. Meanwhile, the domain of the transformed new map becomes narrower at each step in the iterative process. Different from the elimination of mean values by the intersection property, we use the unique property of the normal form. Moreover, in the proof of the convergence of the compositions of infinitely transformations, we take notice of the reversibility of every new mapping, which can be guaranteed by setting the odd or even functions in the right hands of homological equations.

At last, we give an application of the main theorem to show the boundedness of all solutions for the following nonlinear resonant oscillator
 \begin{equation}\label{appl-equa1}
 x''+\varphi(x)f(x')+\omega^2 x+g(x)=p(t),
 \end{equation}
where the functions $f,\varphi,g,p$ satisfy assumptions $(H_1)-(H_5)$ in Section 5. The boundedness conjecture of related problem $ x''+h(t,x)=0$ started from Littlewood \cite{Littwood1}, and in 1976 Morris \cite{Morris} gave a positive answer for $h(t,x)=2x^3+p(t)$. By using KAM theory, he showed the existence of a family of invariant curves and proved the boundedness of every solution of the superlinear case. Additional works for general superlinear and sublinear nonlinearities $h(t,x)$ may be found in \cite{dieckerhoff,Liu,Liu1999,you,kupper,Liu1}. The left class of problems is those with linear growth at infinity, where the case $h(t,x)=n^2x+\varphi(x)-p(t)$ has been settled in \cite{younonlinear,ortega1}. These results are also based on KAM theory.

In this paper, we are concerned with the boundedness of solutions to (\ref{appl-equa1}), where $h(t,x,x')=\varphi(x)f(x')+\omega^2 x+g(x)-p(t)$ has relation with $x'$. Consequently, for $f\neq0$, equation (\ref{appl-equa1}) will no longer be a Hamiltonian system, which means that the classical twist theorem of Moser is not applicable. By setting the parity of $f$ and $p$, we regard (\ref{appl-equa1}) as a reversible system to solve this problem. In \cite{Liubounde}, Kunze, K\"{u}pper and Liu investigated the reversible case $h(t,x,x')=f(x)x'+n^2 x+\varphi(x)-p(t)$ where $p(t)$ is $2\pi$ periodic. By some smoothness hypothesis and a suitable form of small twist theorem, they obtained the boundedness of solutions. A sharp result concerning the unboundedness of solutions was also considered. Later, Li and Ma proved the boundedness of solutions for $h(t,x,x')=f(x')+\omega^2 x+g(x)-p(t)$  with periodic $p(t)$ and $g(+\infty)-g(-\infty)\neq0$ where $g(\pm\infty)=\lim_{x\rightarrow\pm\infty}g(x)$ at $\omega\notin \emph{Q}$ and $\omega\notin \emph{N}$ cases. Combining with the normal form case, we will show similar twist condition $(\varphi(+\infty)-\varphi(-\infty))f(+\infty)+(g(+\infty)-g(-\infty))\neq0$ for (\ref{appl-equa1}) is no longer necessary.

Thus, our second main result is as following.

\begin{theorem}
Suppose that $(H_1)-(H_5)$ hold, and
$\omega$ satisfies
$$|k\omega^{-1}-l|\geq\frac{c_0}{|k|^{\sigma}},$$
where $c_0>0,\sigma>0, k\in\mathbb{Z}\setminus 0, l\in\mathbb{Z}$.
Then for every solution $x(t)$ of (\ref{appl-equa}), we have
$$\sup_{t\in\mathbb{R}}(|x(t)|+|x'(t)|)<+\infty.$$
\end{theorem}

The paper is organized as follows. In Section 2 we recall some notations and give some properties of the normal forms. The iteration lemma and the proof of Theorem 1.1 are given in Sections 3 and 4. As an application, we discuss the boundedness of all solutions for equation \ref{appl-equa} in the last section.

\section{Notations and normal form}

In this section, we will provide some notations and give the definition of normal form without considering the convergence of the power series.

Let $\mathbb{A}:=\mathbb{T}\times\mathbb{R}$ be an infinite cylinder, where $\mathbb{T}=\mathbb{R}/2\pi\mathbb{Z}$. The smooth diffeomorphisms on $\mathbb{A}$ which homotopy to the identity maps are denoted by $Diff_0^\infty(\mathbb{A})$.
We regard as $f(\theta,r)\sim O(r^s)$, if the real analytic function $f(\theta,r)$ can be represented by $$f(\theta,r)=\sum_{k\geq s}f_k(\theta)r^k,$$ and denote  $O(r^s)$ as $O_s$ for short.
We set the total space as
\begin{equation*}
\mathbb{P}=\{\Phi:\mathbb{C}^2\rightarrow\mathbb{C}:\Phi(\theta,r)=\Phi_1(\theta,r)+\Phi_2(\theta,r)+\cdots, \Phi_k\in\mathbb{P}^k,k=1,2,\cdots\},
\end{equation*}
where
\begin{equation*}
\mathbb{P}^k=\{r^k\sum_{j\in\mathbb{Z}^m}\chi_je^{i\langle j,\omega\rangle\theta}:|\Im\theta|<t<1,|r|<\rho<1,\chi_j\in\mathbb{C}, \overline{\chi}_j=\chi_{-j} \},
\end{equation*}
where $\overline{\chi}_j$ is conjugate of $\chi_j$.

Denote the mapping (\ref{main-equa}) in the form of
$$M=I+\beta+F,$$
where $I$ is the identity mapping, $\beta=(\gamma_0,0)^{T}$ and $F=(f(\theta, r),g(\theta, r))^{T}$. Then the linearized normal form (\ref{norm-form}) can be expressed as $I+\beta$.

Introduce an operator $R:\mathbb{P}\rightarrow\mathbb{P}$ as
\begin{equation}\label{r}
R\Phi=\Phi\circ(I+\beta)-\Phi=\Phi(\theta+\gamma_0,r)-\Phi(\theta,r).
\end{equation}
It is obvious that $R$ is a linear mapping and satisfies
$$R|_{\mathbb{P}^k}:\mathbb{P}^k\rightarrow\mathbb{P}^k, k=1,2,\cdots.$$

In the following, we will give the kernel and image of the image $R$.
\begin{lemma}\label{ker-ima}\cite{Li}
Suppose (\ref{dioph}) holds, then the kernel of $R$ is the set of all constant functions in $\theta$, that is
$$\mathbb{K}=\{\Phi:\mathbb{C}^2\rightarrow\mathbb{C}:\Phi=\Phi(r)=\sum_{k=1}^{+\infty}\chi_kr^k,\chi_k\in \mathbb{C}\},$$
and the image of $R$ is the set of all functions which has vanishing average value in $\theta$, that is,
$$\mathbb{M}=\{\Phi:\mathbb{C}^2\rightarrow\mathbb{C}:\Phi=\Phi(\theta,r)=\sum_{k=1}^{+\infty}\sum_{j\in\mathbb{Z}^m\setminus\{0\}}\chi_{kj}e^{i\langle j,\omega\rangle\theta}r^k,\chi_{kj}\in \mathbb{C},\overline{\chi}_{kj}=\chi_{k(-j)}\}.$$
\end{lemma}

From the above discussions,
\begin{equation}\label{dire-sum1}
\mathbb{P}=\mathbb{K}\oplus\mathbb{M},\quad {\mathbb{P}^k}={\mathbb{K}^k}\oplus{\mathbb{M}^k},\quad k=1,2,\cdots,
\end{equation}
and
\begin{equation}\label{comp1}
\quad R|_{\mathbb{P}^k}\circ R^{-1}|_{\mathbb{P}^k}=I|_{\mathbb{M}^k}.
\end{equation}

Similarly, the mapping $R\times R$ on $\mathbb{P}\times\mathbb{P}$ is written in the form
$$(R\times R)(\Phi,\Psi)^{T}=(R\Phi,R\Psi)^{T}=(\Phi\circ(I+\beta)-\Phi,\Psi\circ(I+\beta)-\Psi)^{T}.$$
Moreover, it is easy to derive the properties in (\ref{dire-sum1}) and (\ref{comp1}) for $\mathbb{P}\times\mathbb{P}$ and $R\times R$.

Now we introduce the definition of normal forms and the uniqueness of the normal form, and one can refer to \cite{russman,Liuhu,Li} for details.
\begin{definition}
Assume (\ref{dire-sum1})-(\ref{comp1}) hold. $n$ is called a formal normal form of the mapping $M$, if
\item[{\textrm{(i)}}]: $n-(I+\beta)\in\mathbb{K}\times\mathbb{K}$,

\item[{\textrm{(ii)}}]: there is a change of variables $V=I+\eta$, $\eta\in (R\times R)^{-1}(\mathbb{M}\times\mathbb{M})$, such that $V^{-1}\circ M\circ V=n$.
\end{definition}
\begin{lemma}\label{unique}\cite{Li}
Assume (\ref{dire-sum1})-(\ref{comp1}) hold, if $I+\beta$ is a formal normal form of the given mapping $M$, then the set of formal normal forms of $M$ is $\{I+\beta\}$.
\end{lemma}

In what follows, we will meet the difference equation, the so-called `homological equation':
\begin{equation}\label{homo-1}
Ru(\theta,r)=h(\theta,r),
\end{equation}
where $R$ is defined in (\ref{r}) and $h(\theta,r)$ is a quasi-periodic function in $\mathbb{P}$.
\begin{lemma}\label{lemma-1}
Suppose that $h(\theta,r)$ is real analytic and quasi-periodic in $\theta$ with frequency $\omega=(\omega_1,\omega_2,...,\omega_m)$, the constant $\gamma_0$ satisfies (\ref{dioph}). Then for any $0<\tau<t$,  difference equation (\ref{homo-1}) has unique real analytic quasi-periodic solution $u(\theta,r)$  with frequency $\omega$ and zero average value, that is, $[u]=\lim_{T\rightarrow\infty}\frac{1}{T}\int_0^Tu(\theta,r)d\theta=0$, if and only if $h(\theta,r)\in \mathbb{M}$. Moreover, we have the estimate
$$\|u\|_\tau\leq\frac{c}{(t-\tau)^{m+\sigma}}\|h\|_t,$$
where $\| \cdot\|_a$ denotes the supremum norm in the domain $\{(\theta,r)\in\mathbb{C}^2,|\Im\theta|<a,|r|<1\}$.
If $h(-\theta-\gamma_0,r)=h(\theta,r)$, then $u$ is odd in $\theta$; if $h(-\theta-\gamma_0,r)=-h(\theta,r)$, then $u$ is even in $\theta$.
\end{lemma}

The proof is similar with lemma 2 in \cite{Liu-1}, so we omit it here.
\begin{lemma}\label{prop}
Under the assumptions of diffeomorphism $M$ in (\ref{main-equa}), for every $N\geq2$, there exist a neighborhood of $r=0$ in $\mathbb{A}$, a smooth diffeomorphism $V\in Diff_0^\infty(\mathbb{A})$, and N-1 numbers $\gamma_k\in\mathbb{R}, k=1,2,...,N-1$, such that $\hat{M}=V^{-1}\circ M\circ V$ is in the form of
$$\hat{M}(\theta,r)=\Big(\theta+\gamma_0+\sum_{k=1}^{N-1}\gamma_kr^k+\phi(\theta,r)r^N,r+\varphi(\theta,r)r^N\Big),$$
where $\phi,\varphi$ are smooth mappings. Moreover, $\hat{M}$ is reversible with respect to the involution $G:(\theta,r)\mapsto(-\theta,r)$, and has the invariant curve $r=0$.
\end{lemma}

\noindent{\bf Proof.}  Since $M$ in (\ref{main-equa}) has the invariant curve $r=0$, it can be written as
\begin{equation}
M(\theta,r)=\Big(\theta+\gamma_0+\phi_1(\theta)r,\varphi_1(\theta)r\Big)+O(r^2),
\end{equation}
where $\phi_1,\varphi_1\in C^\infty(\mathbb{T},\mathbb{R})$.
 $M\in Diff_0^\infty(\mathbb{A})$ means that the determinant of Jacobian matrix for $M$ is not vanishing at $r=0$, that is, $\varphi_1(\theta)\neq 0$. Without loss of generality, we may assume $\varphi_1(\theta)>0$.
 In view of (\ref{dioph}) and Lemma \ref{lemma-1},  for every mapping $f\in C^\infty(\mathbb{T}^m,\mathbb{R})$, there is a unique function $g\in C^\infty(\mathbb{T}^m,\mathbb{R})$, such that
\begin{equation*}\label{diff1}
f(\theta)=b+g(\theta)-g(\theta+\gamma_0),\quad\forall\theta\in\mathbb{T}^m,
\end{equation*}
where $b=[f]=\lim_{T\rightarrow\infty}\frac{1}{T}\int_0^T f(\theta)d\theta$.

Hence, for $\varphi_1(\theta)>0$, there exists a $g\in C^\infty(\mathbb{T}^m,\mathbb{R})$, such that,
\begin{equation*}\label{diff11}
\ln{\varphi_1(\theta)}=[\ln{\varphi_1(\theta)}]+g(\theta)-g(\theta+\gamma_0),\quad\forall\theta\in\mathbb{T}^m,
\end{equation*}
that is,
 \begin{equation}\label{diff2}
\varphi_1(\theta)=b_1\frac{g_1^{(2)}(\theta)}{g_1^{(2)}(\theta+\gamma_0)},
\end{equation}
where $b_1=e^{[\ln{\varphi_1(\theta)}]}$, $g_1^{(2)}=e^{g(\theta)}\in C^\infty(\mathbb{T}^m,\mathbb{R})$.

Let
\begin{equation*}
V_1^{(2)}(\theta,r)=\Big(\theta,\frac{r}{\sqrt{g_1^{(2)}(\theta)g_1^{(2)}(-\theta)}} \Big).
\end{equation*}
Then under the transformation $V_1^{(2)}$, the original mapping $M$ is changed into the form
\begin{equation*}
\begin{aligned}
M_1^{(1)}(\theta,r)=&(V_1^{(2)})^{-1}\circ M\circ V_1^{(2)} \\
=&(\theta+\gamma_0+\tilde{\phi_1}(\theta)r,\varphi_1(\theta)\frac{\sqrt{g_1^{(2)}(\theta+\gamma_0)g_1^{(2)}(-\theta-\gamma_0)}}{\sqrt{g_1^{(2)}(\theta)g_1^{(2)}(-\theta)}}\cdot r)+O(r^2),
\end{aligned}
\end{equation*}
where $\tilde{\phi_1}=\phi_1(\theta)/\sqrt{g_1^{(2)}(\theta)g_1^{(2)}(-\theta)}$.

Since $M$ is reversible with respect to the involution $G$, that is,
$$M\circ G\circ M(\theta,r)=G(\theta,r),$$
it is easy to obtain that
$\varphi_1(\theta)\varphi_1(-\theta-\gamma_0)r=r$, then we have $\varphi_1^{-1}(\theta)=\varphi_1(-\theta-\gamma_0)$.

By (\ref{diff2}), one has that
\begin{equation}\label{diff4}
\varphi_1(-\theta-\gamma_0)=b_1\frac{g_1^{(2)}(-\theta-\gamma_0)}{g_1^{(2)}(-\theta)}.
\end{equation}
Hence, combining with (\ref{diff2}) and (\ref{diff4}), we obtain
\begin{equation*}
\varphi_1(\theta)\frac{\sqrt{g_1^{(2)}(\theta+\gamma_0)g_1^{(2)}(-\theta-\gamma_0)}}{\sqrt{g_1^{(2)}(\theta)g_1^{(2)}(-\theta)}}=1.
\end{equation*}
It follows that
\begin{equation*}
M_1^{(1)}(\theta,r)=\Big(\theta+\gamma_0+\tilde{\phi_1}(\theta)r,r\Big)+O(r^2).
\end{equation*}
Finally, the commutativity of $G$ and $V_1^{(2)}$, that is $GV_1^{(2)}=V_1^{(2)}G$, yields that $M_1^{(1)}$ is also reversible with respect to $G$.

From (\ref{dioph}), there is a function $g_1^{(1)}\in C^\infty(\mathbb{T}^m,\mathbb{R})$ satisfying
\begin{equation}\label{diff3}
\tilde{\phi}_1(\theta)=\gamma_1+g_1^{(1)}(\theta)-g_1^{(1)}(\theta+\gamma_0),
\end{equation}
where $\gamma_1=[\tilde{\phi}_1(\theta)]$. We define
\begin{equation*}\label{equa1}
V_1^{(1)}(\theta,r)=(\theta-\frac{g_1^{(1)}(\theta)-g_1^{(1)}(-\theta)}{2}\cdot r,r),
\end{equation*}
then $M_1^{(1)}$ is changed into
\begin{equation}\label{reve-sys1}
\begin{aligned}
M_1(\theta,r)=&(V_1^{(1)})^{-1}\circ M_1^{(1)}\circ V_1^{(1)}\\
=&\Big(\theta+\gamma_0+\big(\tilde{\phi_1}-\frac{g_1^{(1)}(\theta)-g_1^{(1)}(-\theta)}{2}+\frac{g_1^{(1)}(\theta+\gamma_0)-g_1^{(1)}(-\theta-\gamma_0)}{2} \big)\cdot r,r\Big)+O(r^2),
\end{aligned}
\end{equation}

Since $M_1^{(1)}$ is reversible with respect to the involution $G$, that is,
\begin{equation*}\label{rever}
M_1^{(1)}\circ G\circ M_1^{(1)}(\theta,r)=G(\theta,r)=(-\theta,r),
\end{equation*}
the left side of the above equality can be written as
\begin{equation*}
\begin{aligned}
&\Big(-\theta-\tilde{\phi_1}(\theta)r+\tilde{\phi_1}(-\theta-\gamma_0) r,r\Big)+O(r^2),
\end{aligned}
\end{equation*}
which leads to
\begin{equation}\label{equa12}
\tilde{\phi_1}(\theta)=\tilde{\phi_1}(-\theta-\gamma_0).
\end{equation}
A substitution of (\ref{diff3}), (\ref{equa12}) and
$$\tilde{\phi_1}(-\theta-\gamma_0)=\gamma_1+g_1^{(1)}(-\theta-\gamma_0)-g_1^{(1)}(-\theta)$$
into (\ref{reve-sys1}) yields that
\begin{equation}\label{reve-sys2}
M_1(\theta,r)=V_1^{-1}\circ M\circ V_1 =(\theta+\gamma_0+\gamma_1r,r)+O(r^2),
\end{equation}
where $V_1=V_1^{(2)}\circ V_1^{(1)}$.
Due to the commutativity of $G$ and $V_1^{(1)}$, $M_1$ is reversible with respect to $G$.

Next, we expand mapping (\ref{reve-sys2}) to the $r^2$ terms with the form
\begin{equation*}\label{reve-sys3}
M_1(\theta,r)=(\theta+\gamma_0+\gamma_1r+\phi_2(\theta)r^2,r+\varphi_2(\theta)r^2)+O(r^3).
\end{equation*}

Similar to the construction of $V_1^{(1)}$, we can find a constant $b_2$ and a function $g_2^{(2)}\in C^\infty(\mathbb{T}^m,\mathbb{R})$, such that
\begin{equation*}\label{diff5}
\varphi_2(\theta)=b_2+g_2^{(2)}(\theta)-g_2^{(2)}(\theta+\gamma_0).
\end{equation*}
Then under the transformation $$V_2^{(2)}(\theta,r)=\Big(\theta,r-\frac{g_2^{(2)}(\theta)+g_2^{(2)}(-\theta)}{2}\cdot r^2\Big),$$ we have
\begin{equation*}\label{reve-sys4}
\begin{aligned}
&M_2^{(1)}(\theta,r)=(V_2^{(2)})^{-1}\circ M_1\circ V_2^{(2)}\\
&=\Bigg(\theta+\gamma_0+\gamma_1r+\tilde{\phi_2}r^2,r+\Big(\varphi_2(\theta)-\frac{g_2^{(2)}(\theta)+g_2^{(2)}(-\theta)}{2}+\frac{g_2^{(2)}(\theta+\gamma_0)+g_2^{(2)}(-\theta-\gamma_0)}{2}\Big )r^2\Bigg)+O(r^3),
\end{aligned}
\end{equation*}
where $\tilde{\phi_2}(\theta)=\phi_2(\theta)-\gamma_1\frac{g_2^{(2)}(\theta)+g_2^{(2)}(-\theta)}{2}$.
Since $M_1$ is reversible with respect to $G$, we obtain
\begin{equation*}\label{reve-sys5}
M_2^{(1)}(\theta,r)=\Big(\theta+\gamma_0+\gamma_1r+\tilde{\phi_2}(\theta)r^2,r\Big)+O(r^3),
\end{equation*}
and it is also reversible with respect to $G$.

Similarly, there exists a number $\gamma_2$ and a function $g_2^{(1)}\in C^\infty(\mathbb{T}^m,\mathbb{R})$, such that
\begin{equation*}\label{diff6}
\tilde{\phi_2}(\theta)=\gamma_2+g_2^{(1)}(\theta)-g_2^{(1)}(\theta+\gamma_0).
\end{equation*}
By the transformation $$V_2^{(1)}(\theta,r)=\Big(\theta-\frac{g_2^{(1)}(\theta)-g_2^{(1)}(-\theta)}{2}\cdot r^2,r\Big),$$  we have
\begin{equation}\label{reve-sys6}
M_2(\theta,r)=V_2^{-1}\circ M \circ V_2=(\theta+\gamma_0+\gamma_1r+\gamma_2r^2,r)+O(r^3),
\end{equation}
where $V_2=V_1\circ V_2^{(2)}\circ V_2^{(1)} $, and (\ref{reve-sys6}) is also reversible with respect to $G$.

The rest may be deduced by analogy, therefore Lemma \ref{prop} is proved. $\hfill\square$

We denote
\begin{equation*}
\hat{M}(\theta,r)=(\theta+\gamma_0+\sum_{k=1}^{N-1}\gamma_kr^k,r)+O(r^N)
\end{equation*}
for short.

If some $\gamma_k\neq0$, by Moser's twist theorem \cite{moser1}, the mapping $M$ has many invariant curves around the origin if the disturbance terms are sufficiently small. Hence, we only restrict our attention to the existence of invariant curves if $\gamma_k=0,k=1,2,\cdots$.
By the Lemma \ref{prop}, for $N=s_0$ large, there exists a change of variables $V$, transforming (\ref{main-equa})
into $M_{0}=V^{-1}\circ M\circ V$, that is,
\begin{equation}\label{reve-sys71}M_{0}:
\begin{cases}
\theta_1=\theta+\gamma_0+f_{s_0}(\theta,r),\\
r_1=r+g_{s_0}(\theta,r).
\end{cases}
\end{equation}
Moreover, $f_{s_0}(\theta,r)\sim O_{s_0},g_{s_0}(\theta,r)\sim O_{s_0}$, and $M_{0}$ is reversible with respect to the involution $G:(\theta,r)\mapsto(-\theta,r)$. In the rest part, we start with the transformed map  $M_{0}$.

\section{The process of iteration }
In this section, we will give the iteration theorem and its proof, which is of vital importance in the approximation process of $M_{0}$ to its linear part. In the classical iteration method, we use the twist condition and intersection or area preserving properties to eliminate the mean values generated by perturbations in the direction of angular and action variables. But for reversible mapping $M_{0}$, there is no twist condition (all $\gamma_k=0$ ), and no intersection or area preserving properties. How do we eliminate the mean values? Fortunately, by R\"{u}ssmann's method in \cite{russman}, we use the uniqueness of the formal normal forms to achieve this goal. Under the conditions of Theorem \ref{main-theo}, assume there is a formal change of variables transforming $M_{0}$ into the linearized form $I+\beta$.

\begin{lemma}\label{exist}
Assume the conditions of Theorem \ref{main-theo} hold, then we can derive a sequence of changes of variables $\{T_n\}$, such that
\begin{equation}\label{seque}
T_n^{-1}\circ M_{0}\circ T_n-(I+\beta)=O_{s_n},
\end{equation}
where $s_n=2^{\alpha+n}+1$, $\alpha $ is a sufficiently large positive integer to be determined later and $T_0=I$. Moreover, for every $n$, the transformed mapping $M_n:=T_n^{-1}\circ M_{0}\circ T_n$ is reversible with respect to the involution $G:(\theta,r)\mapsto(-\theta,r)$.
\end{lemma}

\noindent{\bf Proof.} For $n=0$, we choose $s_0=2^\alpha+1$. By (\ref{reve-sys71}), it is easy to see (\ref{seque}) holds.

Suppose for all $n$, there exists $T_{n}$, such that
$$T_n^{-1}\circ M_{0}\circ T_n-(I+\beta)\sim O_{s_n},$$
then for $(n+1)-$th step, we need to find a change of variables $T_{n+1}$, such that
$$T_{n+1}^{-1}\circ M_{0}\circ T_{n+1}-(I+\beta)\sim O_{s_{n+1}}.$$
Denote $T_n^{-1}\circ M_{0}\circ T_n=M_n$, and $T_{n+1}=T_{n}\circ \Delta T_{n+1}$.

It is in a position to establish a transformation $\Delta T_{n+1}$, satisfying
$$\Delta T_{n+1}^{-1}\circ M_n\circ \Delta T_{n+1}=M_{n+1}.$$
In the sequel, we denote $s_n$ and $\Delta T_{n+1}$ by $s$ and $ \Delta T$, respectively.

Assume the change of variables $\Delta T$ has the form
\begin{equation*}\label{form-chan1}\Delta T:
\begin{cases}
\theta=\xi+u(\xi,\eta),\\
r=\eta+v(\xi,\eta),
\end{cases}
\end{equation*}
 and we denote $M_{n+1}$ as
     \begin{equation*}\label{reve-sys9}M_{n+1}:
\begin{cases}
\xi_1=\xi+\gamma_0+f_{\hat{s}}(\theta,r),\\
\eta_1=\eta+g_{\hat{s}}(\theta,r).
\end{cases}
\end{equation*}
From $M_{n+1}=\Delta T^{-1}\circ M_n\circ \Delta T$, we have $M_n\circ \Delta T=\Delta T\circ M_{n+1}$, it follows that
\begin{equation}\label{equa1}
\begin{cases}
f_{\hat{s}}(\theta,r)=u(\xi,\eta)-u(\xi+\gamma_0+f_{\hat{s}},\eta+g_{\hat{s}})+f_s(\xi+u,\eta+v),\\
g_{\hat{s}}(\theta,r)=v(\xi,\eta)-v(\xi+\gamma_0+f_{\hat{s}},\eta+g_{\hat{s}})+g_s(\xi+u,\eta+v),
\end{cases}
\end{equation}
that is,
\begin{equation}\label{equa2}
\begin{cases}
f_{\hat{s}}(\theta,r)=-(u(\xi+\gamma_0,\eta)-u(\xi,\eta))+u(\xi+\gamma_0,\eta)-u(\xi+\gamma_0+f_{\hat{s}},\eta+g_{\hat{s}})+f_s(\xi+u,\eta+v),\\
g_{\hat{s}}(\theta,r)=-(v(\xi+\gamma_0,\eta)-v(\xi,\eta))+v(\xi+\gamma_0,\eta)-v(\xi+\gamma_0+f_{\hat{s}},\eta+g_{\hat{s}})+g_s(\xi+u,\eta+v).
\end{cases}
\end{equation}

The important point to note here is that we need to ensure $M_{n+1}$ is also reversible under transformation $\Delta T$. According to the concept of reversible mappings, it implies the change of variable $\Delta T$  need to commutes with the involution $G(\xi,\eta)=(-\xi,\eta)$. This requires functions $u$ and $v$ satisfying
\begin{equation}\label{odd-even}
u(-\xi,\eta)=-u(\xi,\eta),\quad\quad v(-\xi,\eta)=v(\xi,\eta).
\end{equation}

From Lemma \ref{lemma-1}, we try to decide the functions $u$ and $v$ from the following modified difference equation:
\begin{equation}\label{equa22}
\begin{cases}
u(\xi+\gamma_0,\eta)-u(\xi,\eta)=\{p_s\}_{\mathbb{M}},\\
 v(\xi+\gamma_0,\eta)-v(\xi,\eta)=\{q_s\}_{\mathbb{M}},
 \end{cases}
\end{equation}
where
\begin{equation}\label{equa23}
p_s(\xi,\eta)=\frac{1}{2}(f_s(\xi,\eta)+f_s(-\xi-\gamma_0,\eta))\quad q_s(\xi,\eta)=\frac{1}{2}(g_s(\xi,\eta)-g_s(-\xi-\gamma_0,\eta)).
\end{equation}
It is easy to verify that $p_s(-\xi-\gamma_0,\eta)=p_s(\xi,\eta),q_s(-\xi-\gamma_0,\eta)=-p_s(\xi,\eta)$.  $\{p_s\}_{\mathbb{M}}$ and $\{q_s\}_{\mathbb{M}}$ also possess this property and $[\{p_s\}_{\mathbb{M}}](\eta)=[\{q_s\}_{\mathbb{M}}](\eta)=0$. So the functions $u$ and $v$
meet the condition (\ref{odd-even}).

Thus, (\ref{equa2}) can be rewritten as
\begin{equation}\label{equa2-3}
\begin{cases}
f_{\hat{s}}(\theta,r)=-\{p_s\}_{\mathbb{M}}+u(\xi+\gamma_0,\eta)-u(\xi+\gamma_0+f_{\hat{s}},\eta+g_{\hat{s}})+f_s(\xi+u,\eta+v)+p_s(\xi,\eta)-p_s(\xi,\eta),\\
g_{\hat{s}}(\theta,r)=-\{q_s\}_{\mathbb{M}}+v(\xi+\gamma_0,\eta)-v(\xi+\gamma_0+f_{\hat{s}},\eta+g_{\hat{s}})+g_s(\xi+u,\eta+v)+q_s(\xi,\eta)-q_s(\xi,\eta).
\end{cases}
\end{equation}

Since $f_s\sim O_s,g_s\sim O_s$, it follows that $p_s\sim O_s, q_s\sim O_s$, and by (\ref{equa1}) and (\ref{equa22}), there are $$u\sim O_s, v\sim O_s,f_{\hat{s}}\sim O_s , g_{\hat{s}}\sim O_s,$$
and
\begin{equation*}
\begin{aligned}
&u(\xi+\gamma_0+f_{\hat{s}},\eta+g_{\hat{s}})-u(\xi+\gamma_0,\eta)\\
=&D_{\xi}u(\xi+\gamma_0,\eta)\cdot f_{\hat{s}}+D_{\eta}u(\xi+\gamma_0,\eta)\cdot g_{\hat{s}}\\
\sim& O_{2s-1}.
\end{aligned}
\end{equation*}

 In the same way, we have$$v(\xi+\gamma_0+f_{\hat{s}},\eta+g_{\hat{s}})-v(\xi+\gamma_0,\eta)\sim O_{2s-1},$$
 $$f_s(\xi+u,\eta+v)-f_s(\xi,\eta)\sim O_{2s-1},$$
 $$g_s(\xi+u,\eta+v)-g_s(\xi,\eta)\sim O_{2s-1},$$
 $$f_s(-\xi-\gamma_0-f_s,\eta+g_s)-f_s(-\xi-\gamma_0,\eta)\sim  O_{2s-1}.$$
In the following, we will prove
\begin{equation}\label{equa24}
f_s(\xi+u,\eta+v)-p_s(\xi,\eta)\sim O_{2s-1},\quad g_s(\xi+u,\eta+v)-q_s(\xi,\eta)\sim O_{2s-1}.
\end{equation}
From the definitions of $p_s$ and $q_s$, we have
$$f_s(\xi+u,\eta+v)-p_s(\xi,\eta)
=\frac{1}{2}(f_s(\xi+u,\eta+v)-f_s(\xi,\eta))+\frac{1}{2}(f_s(\xi+u,\eta+v)-f_s(-\xi-\gamma_0,\eta)),$$
$$g_s(\xi+u,\eta+v)-q_s(\xi,\eta)\\
=\frac{1}{2}(g_s(\xi+u,\eta+v)-g_s(\xi,\eta))+\frac{1}{2}(g_s(\xi+u,\eta+v)+g_s(-\xi-\gamma_0,\eta)).$$

Since the mapping  $M_n$ is reversible, i.e., $M_nGM_n=G$, it follows that
\begin{equation}\label{reve-prop1}
f_s(-\xi-\gamma_0-f_s,\eta+g_s)-f_s(\xi,\eta)=0,\ \ \
g_s(-\xi-\gamma_0-f_s,\eta+g_s)+g_s(\xi,\eta)=0.
\end{equation}
Hence, by (\ref{reve-prop1}), we have
\begin{equation*}
\begin{aligned}&f_s(\xi+u,\eta+v)-f_s(-\xi-\gamma_0,\eta)\\
&=f_s(\xi+u,\eta+v)-f_s(-\xi-\gamma_0-f_s,\eta+g_s)+f_s(-\xi-\gamma_0-f_s,\eta+g_s)-f_s(-\xi-\gamma_0,\eta)\\
&=f_s(\xi+u,\eta+v)-f_s(\xi,\eta)+f_s(-\xi-\gamma_0-f_s,\eta+g_s)-f_s(-\xi-\gamma_0,\eta)\\
&\sim  O_{2s-1}.
\end{aligned}
\end{equation*}
By the same reason, $g_s(\xi+u,\eta+v)-g_s(-\xi-\gamma_0,\eta)\sim O_{2s-1}$.
As a consequence, (\ref{equa24}) holds.

By (\ref{equa2-3}), we have
\begin{equation*}\label{equa7}
\begin{cases}
f_{\hat{s}}(\xi,\eta)=\{p_s\}_{\mathbb {K}}+O_{2s-1},\\
g_{\hat{s}}(\xi,\eta)=\{q_s\}_{\mathbb {K}}+O_{2s-1}.
\end{cases}
\end{equation*}
Then it follows that
\begin{equation}\label{equa8}
M_{n+1}=\Delta T^{-1}\circ M_n\circ \Delta T=I+\beta+\{H\}_{\mathbb {K}\times\mathbb{K}}+O_{2s-1},
\end{equation}
where $H=(p_s,q_s),\hat{s}=2s-1$.

For (\ref{equa22}) and (\ref{equa8}), nothing changes if we replace $H$ by
$$H^*=(p_s^*,q_s^*)=\sum_{k=s}^{2s-2}H_k,$$
 that is,
 \begin{equation*}\label{equa88}
(u(\xi,\eta),v(\xi,\eta))=R^{-1}\{H^*\}_{\mathbb{M}\times\mathbb{M}},
\end{equation*}
where $Ru=u(\xi+\gamma_0,\eta)-u(\xi,\eta)$,
instead of (\ref{equa22}), and
\begin{equation}\label{equa9}
M_{2s-1}=\Delta T^{-1}\circ M_s\circ \Delta T=I+\beta+\{H^*\}_{\mathbb {K}\times\mathbb{K}}+O_{2s-1},
\end{equation}
instead of (\ref{equa8}).

In the following, we  show that $\{H^*\}_{\mathbb {K}\times\mathbb{K}}=0$.
Otherwise, if $\{H^*\}_{\mathbb {K}\times\mathbb{K}}=K_m+O_{m+1}$, $K_m\in K^{m}\times K^{m}, K_m\neq0, s\leq m\leq 2s-2$.  Then from (\ref{equa9}), we have
 $$\Delta T^{-1}\circ M_s\circ \Delta T=I+\beta+K_m+O_{m+1}+O_{2s-1},$$
which contradicts our assumption for $M_0$ and Lemma \ref{unique}.
Hence, $\{H^*\}_{\mathbb {K}\times\mathbb{K}}=0$.
As a consequence, $$f_{2s-1}(\xi,\eta)\sim O_{2s-1},
\quad g_{2s-1}(\xi,\eta)\sim O_{2s-1},$$
and there exists a $T_{n+1}=T_{n}\circ \Delta T_{n+1}$ such that (\ref{seque}) holds. Therefore the Lemma is proved.$\hfill\square$

Now, we introduce the iteration lemma, which is used infinitely times to transforming $M_0$ close to $I+\beta$. Since the iteration lemma is one step in the iteration process, we write $s$ instead of $s_n$ and $\widetilde{M}$ instead of mapping $M_n$ in (\ref{seque}).

Set several complex domains
$$D=\{(\theta,r): |\Im\theta|<t,|r|<\rho\},$$
$$B=\{(\theta,r): |\Im\theta|<\tau,|r|<\varrho\},$$
$$B^{(k)}=\{(\theta,r): |\Im\theta|<t-\frac{k(t-\tau)}{4},|r|<\rho-\frac{k(\rho-\varrho)}{4}\},k=1,2,3,$$
with $\tau<t<1,\varrho<\rho<1$.
It is easy to get $B\subset B^{(3)}\subset B^{(2)}\subset B^{(1)}\subset D$.
In what follows, we denote the norm $|f|_D=\sup_{(\theta,r)\in D}|f(\theta,r)|$.

\begin{lemma}(Iteration Lemma)\label{iter}
Suppose that the real analytic mapping
\begin{equation}\label{reve-sys7}\widetilde{M}:
\begin{cases}
\theta_1=\theta+\gamma_0+f_{s}(\theta,r),\\
r_1=r+g_{s}(\theta,r),
\end{cases}
\end{equation}
where $$f_{s}(\theta,r)=\sum_{k\geq s}^{\infty}\sum_{j\in \mathbb{Z}^l}f_{kj}e^{i\langle j,\omega\rangle\theta}r^k, \quad g_{s}(\theta,r)=\sum_{k\geq s}^{\infty}\sum_{j\in \mathbb{Z}^l}g_{kj}e^{i\langle j,\omega\rangle\theta}r^k,$$
satisfy

\item[{\textrm{(i)}}]: $\widetilde{M}$ is reversible with respect to $G:(\theta,r)\rightarrow(-\theta,r)$, and the constant $\gamma_0$ satisfies (\ref{dioph});

\item[{\textrm{(ii)}}]: $|f_s|_D+|g_s|_D<d, d<\min\{\frac{t-\tau}{4},\frac{\rho-\varrho}{4}\}$;

\item[{\textrm{(iii)}}]: $\nu=c_7d(t-\tau)^{-m-\sigma}(\rho-\varrho)^{-1}(\frac{1}{t-\tau}+\frac{1}{\rho-\varrho})<\min\{\frac{t-\tau}{4},\frac{\rho-\varrho}{4},\frac{1}{5}\}$,
    where $c_7$ and $c_i(i=1,2,3,4,5,6,8,9)$ are  positive constants, which will determined later, and depend only on $c_0, \omega, \sigma$.

     Then there exists a change of variables $\Delta T$
 such that $\widetilde{M}$ is transformed into $\overline{M}=\Delta T^{-1}\circ \widetilde{M}\circ \Delta T$ with the form
     \begin{equation}\label{reve-sys8}\overline{M}:
\begin{cases}
\xi_1=\xi+\gamma_0+f_{2s-1}(\xi,\eta),\\
\eta_1=\eta+g_{2s-1}(\xi,\eta),
\end{cases}
\end{equation}
and $\overline{M}$ is reversible with respect to $G$. Moreover,
the following estimates hold:
\begin{equation}\label{esti3}
|u|_{B^{(1)}}+|v|_{B^{(1)}}<c_2d(t-\tau)^{-m-\sigma}(\rho-\varrho)^{-1},
\end{equation}
\begin{equation}\label{esti4}
|u_\xi|_{B^{(1)}}+|u_\eta|_{B^{(1)}}+|v_\xi|_{B^{(1)}}+|v_\eta|_{B^{(1)}}<c_3d(t-\tau)^{-m-\sigma}(\rho-\varrho)^{-1}(\frac{1}{t-\tau}+\frac{1}{\rho-\varrho}),
\end{equation}
\begin{equation}\label{esti5}
|f_{2s-1}(\xi,\eta)|_{B}+|g_{2s-1}(\xi,\eta)|_{B}\leq\frac{c_6(\rho-\varrho)^{-1}(de^{\frac{1}{2}(-2s-1)(\rho-\varrho)}+d^2(t-\tau)^{-m-\sigma}(\frac{1}{t-\tau}+\frac{1}{\rho-\varrho}))}{1-\nu}.
\end{equation}
\end{lemma}
\noindent{\bf Proof.} Firstly, the existence of such a change
\begin{equation}\label{form-chan1}\Delta T:
\begin{cases}
\theta=\xi+u(\xi,\eta),\\
r=\eta+v(\xi,\eta)
\end{cases}
\end{equation}
is guaranteed by Lemma \ref{exist}. Secondly, we give some estimates of $\Delta T$ and $\overline{M}$.

Suppose $u,v$ in (\ref{form-chan1}) and $p_s$, $q_s$ in (\ref{equa23}) have expansions of the type:
\begin{equation}\label{equa91}
u(\xi,\eta)=\sum_{k=s}^{2s-2}\sum_{j\in\mathbb{Z}^m}u_{kj}e^{i\langle j,\omega\rangle\xi}\eta^k,\quad v(\xi,\eta)=\sum_{k=s}^{2s-2}\sum_{j\in\mathbb{Z}^m}v_{kj}e^{i\langle j,\omega\rangle\xi}\eta^k,
\end{equation}
$$p_s(\xi,\eta)=\sum_{k=s}^{+\infty}\sum_{j\in\mathbb{Z}^m}p_{kj}e^{i\langle j,\omega\rangle\xi}\eta^k,\quad q_s(\xi,\eta)=\sum_{k=s}^{+\infty}\sum_{j\in\mathbb{Z}^m}q_{kj}e^{i\langle j,\omega\rangle\xi}\eta^k.$$
Then, \begin{equation}\label{equa92}
\{p_s^*\}_{\mathbb{M}}=\sum_{k=s}^{2s-2}\sum_{j\in\mathbb{Z}^m\setminus\{0\}}p_{kj}e^{i\langle j,\omega\rangle\xi}\eta^k,\quad \{q_s^*\}_{\mathbb{M}}=\sum_{k=s}^{2s-2}\sum_{j\in\mathbb{Z}^m\setminus\{0\}}q_{kj}e^{i\langle j,\omega\rangle\xi}\eta^k.\end{equation}
Substituting (\ref{equa91}) and (\ref{equa92}) into (\ref{equa22}), we have
$$u_{kj}(e^{i\langle j,\omega\rangle\gamma_0}-1)=p_{kj}, \quad v_{kj}(e^{i\langle j,\omega\rangle\gamma_0}-1)=q_{kj},\quad s\leq k\leq 2s-1,\ j\in\mathbb{Z}^m\setminus\{0\},$$
hence,
 \begin{equation*}\label{equa93}
 u_{kj}=\frac{p_{kj}}{e^{i\langle j,\omega\rangle\gamma_0}-1},\quad v_{kj}=\frac{q_{kj}}{e^{i\langle j,\omega\rangle\gamma_0}-1},\quad s\leq k\leq 2s-1, \ j\in\mathbb{Z}^m\setminus\{0\},
\end{equation*}
and $u_{k0}=0,v_{k0}=0$.

Since the constant $\gamma_0$ and $\omega$ satisfy (\ref{dioph}),  it means that
 \begin{equation*}\label{esti11}
 |e^{i\langle j,\omega\rangle\gamma_0}-1|\geq\frac{4c_0}{|j|^{\sigma}}.
\end{equation*}
Denote $p_s=\sum_{k=s}^{+\infty}p_k(\xi)\eta^k$, $p_k(\xi)=\sum_{j\in\mathbb{Z}^m}p_{kj}e^{i\langle j,\omega\rangle\xi}$, then
$$p_s^*=\sum_{k=s}^{2s-2}p_k(\xi)\eta^k,\quad p_k(\xi)=\sum_{j\in\mathbb{Z}^m}p_{kj}e^{i\langle j,\omega\rangle\xi}.$$
By Cauchy's estimate and the analyticity of $p_s$ and $q_s$, we have
 $$|p_k|_D\leq d\rho^{-k},\ \ |p_{kj}|_D\leq d\rho^{-k}e^{-|j||\omega|t}.$$
Denote a narrower strip $$D^*=\{(\xi,\eta):|\Im\xi|<t-\delta_1,,|\eta|<\rho-\delta_2\},$$where
$\delta_1=\frac{t-\tau}{5},\delta_2=\frac{\rho-\varrho}{5}$, and it is easy to prove $B\subset B^{(3)}\subset B^{(2)} \subset
B^{(1)}\subset D^*\subset D$.

Then we have the estimate
 \begin{equation*}\label{esti12}
 \begin{aligned}
 |u|_{D^*}=&\left|\sum_{k=s}^{2s-2}\sum_{j\in\mathbb{Z}^m\setminus\{0\}}\frac{p_{kj}}{e^{i\langle j,\omega\rangle\gamma_0}-1}e^{i\langle j,\omega\rangle\xi}\eta^k\right|\\
\leq&\sum_{k=s}^{2s-2}\sum_{j\in\mathbb{Z}^m\setminus\{0\}}\frac{|j|^\sigma}{4c_0}d\rho^{-k}e^{-|j||\omega|t}e^{|j||\omega|(t-\delta_1)}(\rho-\delta_2)^k\\
=&\frac{d}{4c_0}\sum_{k=s}^{2s-2}(1-\frac{\delta_2}{\rho})^k\sum_{j\in\mathbb{Z}^m\setminus\{0\}}|j|^{\sigma}e^{-|j||\omega|\delta_1}\\
<&c_1d\delta_1^{-m-\sigma}\frac{\rho}{\delta_2}<\frac{1}{2}c_2d(t-\tau)^{-m-\sigma}(\rho-\varrho)^{-1},
 \end{aligned}
\end{equation*}
where $c_2>10c_1$, and $c_1,c_2$ are positive constants depending on $c_0, \sigma,\omega$.
In a similar way,
 \begin{equation*}\label{esti13}
  |v|_{D^*}<\frac{1}{2}c_2d(t-\tau)^{-m-\sigma}(\rho-\varrho)^{-1}.
\end{equation*}
From the above discussions, we have (\ref{esti3}), and
by Cauchy's estimate, we get (\ref{esti4}):
$$|u_\xi|_{D^*}+|u_\eta|_{D^*}+|v_\xi|_{D^*}+|v_\eta|_{D^*}<c_3d(t-\tau)^{-m-\sigma}(\rho-\varrho)^{-1}(\frac{1}{t-\tau}+\frac{1}{\rho-\varrho}),$$
where $c_3>c_2$ is a positive constant.

The last step is to estimate $f_{2s-1}$ and $g_{2s-1}$, which satisfy the equation (\ref{equa2}), that is
\begin{equation*}
\begin{cases}
f_{2s-1}(\xi,\eta)=-\{p_s^*\}_{\mathbb{M}}+u(\xi+\gamma_0,\eta)-u(\xi+\gamma_0+f_{2s-1},\eta+g_{2s-1})+f_s(\xi+u,\eta+v),\\
g_{2s-1}(\xi,\eta)=-\{q_s^*\}_{\mathbb{M}}+v(\xi+\gamma_0,\eta)-v(\xi+\gamma_0+f_{2s-1},\eta+g_{2s-1})+g_s(\xi+u,\eta+v).
\end{cases}
\end{equation*}

Since \begin{equation*}
 \begin{aligned}
 &f_s(\xi+u,\eta+v)-\{p_s^*\}_{\mathbb{M}}\\
 &=\big(f_s(\xi+u,\eta+v)-f_s(\xi,\eta)\big)+\big(f_s(\xi,\eta)-p_s(\xi,\eta)\big)+\big(p_s(\xi,\eta)-p_s^*(\xi,\eta)\big)+\big(p_s^*(\xi,\eta)-\{p_s^*\}_{\mathbb{M}}\big),
 \end{aligned}
\end{equation*}
we divide the estimate of $|f_s(\xi+u,\eta+v)-\{p_s^*\}_{\mathbb{M}}|_{B}$ into four parts.

By (\ref{esti3}), it follows that
\begin{equation*}
 \begin{aligned}
 &|f_s(\xi+u,\eta+v)-f_s(\xi,\eta)|_{B}\\
 &\leq|D_\xi f_s|\cdot |u|+|D_\eta f_s|\cdot |v|\\
 &\leq c_4d^2(t-\tau)^{-m-\sigma}(\rho-\varrho)^{-1}(\frac{1}{t-\tau}+\frac{1}{\rho-\varrho}),
  \end{aligned}
\end{equation*}
where $c_4>c_3$ is a positive constant.
Combining with  (\ref{equa23}) and (\ref{reve-prop1}), we have
\begin{equation*}
 \begin{aligned}
 &|f_s(\xi,\eta)-p_s(\xi,\eta)|_{B}\\
 =&\frac{1}{2}|f_s(\xi,\eta)-f_s(-\xi-\gamma_0,\eta)|\\
  =&\frac{1}{2}|f_s(-\xi-\gamma_0-f_s,\eta+g_s)-f_s(-\xi-\gamma_0,\eta)|\\
 \leq& c_4d^2(\frac{1}{t-\tau}+\frac{1}{\rho-\varrho}).
  \end{aligned}
\end{equation*}
In view of
$p_s=\sum_{k=s}^{+\infty}p_k(\xi)\eta^k$
and $$|p_k|_{B^{(2)}}\leq |p_s|_D\left(\frac{\rho+\varrho}{2}\right)^{-k}<d(\frac{\rho+\varrho}{2})^{-k},$$
it yields
\begin{equation*}
 \begin{aligned}
 &|p_s(\xi,\eta)-p_s^*(\xi,\eta)|_{B}\\
 =&\left|\sum_{k=2s-1}^{+\infty}p_{k}(\xi)\eta^k\right|\leq\sum_{k=2s-1}^{+\infty}|p_k|_{B^{(2)}}|\varrho|^k\\
 \leq&\sum_{k=2s-1}^{+\infty}d\left(\frac{\rho+\varrho}{2}\right)^{-k}|\varrho|^k\\
 =&d\left(1-\frac{\rho-\varrho}{\rho+\varrho}\right)^{2s-1}\sum_{k=0}^{+\infty}(1-\frac{\rho-\varrho}{\rho+\varrho})^{k}\\
 \leq&c_5de^{-\frac{1}{2}(2s-1)(\rho-\varrho)}(\rho-\varrho)^{-1},
  \end{aligned}
\end{equation*}
where $c_5>c_4$ is a positive constant.

From these estimates, we have
$$|f_s(\xi+u,\eta+v)-\{p_s^*\}_{\mathbb{M}}|_{B}<\frac{1}{2}c_6\Big(\rho-\varrho)^{-1}(de^{\frac{1}{2}(-2s-1)(\rho-\varrho)}+d^2(t-\tau)^{-m-\sigma}(\frac{1}{t-\tau}+\frac{1}{\rho-\varrho})\Big),$$
where $c_6>2c_5$ is a positive constant.

As a consequence,
\begin{equation*}
 \begin{aligned}
&|f_{2s-1}(\xi,\eta)|_{B}\\
\leq&\frac{1}{2}c_6(\rho-\varrho)^{-1}(de^{\frac{1}{2}(-2s-1)(\rho-\varrho)}+d^2(t-\tau)^{-m-\sigma}(\frac{1}{t-\tau}+\frac{1}{\rho-\varrho}))\\
&+|u(\xi+\gamma_0,\eta)-u(\xi+\gamma_0+f_{2s-1},\eta+g_{2s-1})|\\
\leq&\frac{1}{2}c_6(\rho-\varrho)^{-1}(de^{\frac{1}{2}(-2s-1)(\rho-\varrho)}+d^2(t-\tau)^{-m-\sigma}(\frac{1}{t-\tau}+\frac{1}{\rho-\varrho}))\\
&+\frac{1}{2}c_7d(t-\tau)^{-m-\sigma}(\rho-\varrho)^{-1}(\frac{1}{t-\tau}+\frac{1}{\rho-\varrho})(|f_{2s-1}|_{B}+|g_{2s-1}|_{B}),
  \end{aligned}
\end{equation*}
where $c_7$ is a positive constant. Similarly, there is same estimate for $|g_{2s-1}(\xi,\eta)|_{B}$.

From the above discussion, we obtain
$$|f_{2s-1}(\xi,\eta)|_{B}+|g_{2s-1}(\xi,\eta)|_{B}\leq\frac{c_6(\rho-\varrho)^{-1}(de^{\frac{1}{2}(-2s-1)(\rho-\varrho)}+d^2(t-\tau)^{-m-\sigma}(\frac{1}{t-\tau}+\frac{1}{\rho-\varrho}))}{1-\nu},$$
by denoting $\nu=c_7d(t-\tau)^{-m-\sigma}(\rho-\varrho)^{-1}(\frac{1}{t-\tau}+\frac{1}{\rho-\varrho})$. The proof is completed.$\hfill\square$

Obviously, applying the iteration lemma to $M_n$ in (\ref{seque}), we obtain $M_{n+1}$ with estimates (\ref{esti3})-(\ref{esti5}), so that iteration process can continue. The specific steps and the convergence of composed mappings are shown in the next section.
\section{Proof of Theorem \ref{main-theo}}

Set some sequences of variables and domains:
$$t_n=\frac{t_0}{2}\big(1+(\frac{2}{3})^n\big),t_0<1,t=t_n,\tau=t_{n+1},$$
$$\rho_n=\frac{t_0}{2}\big(1+(\frac{2}{3})^n\big),\rho_0<1,\rho=\rho_n,\varrho=\rho_{n+1},$$
$$d_{n+1}=(\frac{3}{2})^nd_n^{\frac{4}{3}},d=d_0<1,$$
$$e_n=(\frac{3}{2})^{3n+9}d_n, e_{n+1}=e_n^{\frac{4}{3}},$$
$$D_n=\{(\xi,\eta):|\Im\xi|<t_n,|\eta|<\rho_n\}.$$

In this section, we will verify that there exists a convergent change of variables, transforming (\ref{main-equa}) into (\ref{form-chan}). For this purpose, we need to prove for every $n$, there is a transformation $T_n$ such that $M_0$ is transformed into $M_n$,
     \begin{equation*}\label{reve-sys81}M_{n}:
\begin{cases}
\xi_1=\xi+\gamma_0+f_{s_n}(\xi,\eta),\\
\eta_1=\eta+g_{s_n}(\xi,\eta),
\end{cases}
\end{equation*}
and
\begin{equation}\label{pert-esti}
|f_{s_n}|_{D_n}+|g_{s_n}|_{D_n}<d_n.
\end{equation}
By Lemma \ref{iter}, the existence of transformation $T_n$ is obtained. Thus we just check (\ref{pert-esti}) for all $n$.

When $n=0$, $|f_{s_0}|_{D_0}+|g_{s_0}|_{D_0}<d=d_0$.
Supposing for all $n$, the nonlinear part of $M_n$ satisfies (\ref{pert-esti}), it is in a position to prove that
$$|f_{s_{n+1}}|_{D_{n+1}}+|g_{s_{n+1}}|_{D_{n+1}}<d_{n+1}.$$
Firstly, we have to guarantee the conditions of Lemma \ref{iter}, i.e.,
$$d_n<\min\{\frac{t_n-t_{n+1}}{4},\frac{\rho_n-\rho_{n+1}}{4}\},$$
$$\nu=c_7d_n(t_n-t_{n+1})^{-m-\sigma}(\rho_n-\rho_{n+1})^{-1}(\frac{1}{t_n-t_{n+1}}+\frac{1}{\rho_n-\rho_{n+1}})<\min\{\frac{t_n-t_{n+1}}{4},\frac{\rho_n-\rho_{n+1}}{4},\frac{1}{5}\},$$
which means that
$$(\frac{3}{2})^nd_n<\min\{\frac{t_0}{24},\frac{\rho_0}{24}\},$$
and
$$c_7d_n(\frac{t_0}{6})^{-m-\sigma}(\frac{\rho_0}{6})^{-1}\Big((\frac{t_0}{6})^{-1}+(\frac{\rho_0}{6})^{-1}\Big)(\frac{3}{2})^{n(m+\sigma+2)}<\min\{\frac{t_0}{24},\frac{\rho_0}{24},\frac{1}{5}\}.$$
From the definition of $d_n$, we have $d_n=(\frac{2}{3})^{3n+9}\Big((\frac{3}{2})^9d_0\Big)^{(\frac{4}{3})^n}$.  Obviously, $d_0$ depends on $t_0,\rho_0,m,\sigma$ and $c_0$, thus we can choose $\alpha$ large enough such that $d_0$ is sufficiently small, finally all above inequalities hold.
By Lemma \ref{iter}, we derive
\begin{equation*}
\begin{aligned}&|f_{s_{n+1}}(\xi,\eta)|_{D_{n+1}}+|g_{s_{n+1}}
(\xi,\eta)|_{D_{n+1}}\\
&\leq\frac{c_6(\rho_n-\rho_{n+1})^{-1}(d_ne^{\frac{1}{2}(2s_n-1)(\rho_n-\rho_{n+1})}+(d_n)^2(t_n-t_{n+1})^{-m-\sigma}(\frac{1}{t_n-t_{n+1}}+\frac{1}{\rho_n-\rho_{n+1}}))}{1-\nu}.
\end{aligned}
\end{equation*}

On the one hand, we will prove the inequality
$$\frac{c_6(\rho_n-\rho_{n+1})^{-1}d_ne^{-\frac{1}{2}(2s_n-1)(\rho_n-\rho_{n+1})}}{1-\nu}<\frac{1}{2}d_{n+1}.$$
Since
\begin{equation*}
\begin{aligned}&\frac{c_6(\rho_n-\rho_{n+1})^{-1}d_ne^{-\frac{1}{2}(2s_n-1)(\rho_n-\rho_{n+1})}}{1-\nu}\\
&<\frac{5}{4}c_6(\frac{\rho_0}{6})^{-1}(\frac{3}{2})^ne^{-\frac{\rho_0}{12}(\frac{2}{3})^n(2s_n-1)}d_n\\
&=\frac{5}{4}c_6(\frac{\rho_0}{6})^{-1}e^{-\frac{\rho_0}{12}(\frac{2}{3})^n(2s_n-1)}d_n^{-\frac{1}{3}}(\frac{3}{2})^nd_n^{\frac{4}{3}}\\
&<\frac{1}{2}\Big(c_8(\frac{\rho_0}{6})^{-1}e^{-\frac{\rho_0}{12}(\frac{2}{3})^n(2s_n-1)}d_n^{-\frac{1}{3}}\Big)d_{n+1},
\end{aligned}
\end{equation*}
where $c_8>\frac{5}{2}c_6$ is a positive constant,
we will only verify the term $c_8(\frac{\rho_0}{6})^{-1}e^{-\frac{\rho_0}{12}(\frac{2}{3})^n(2s_n-1)}d_n^{-\frac{1}{3}}<1$.
In fact
\begin{equation*}
\begin{aligned}
&c_8(\frac{\rho_0}{6})^{-1}e^{-\frac{\rho_0}{12}(\frac{2}{3})^n(2s_n-1)}d_n^{-\frac{1}{3}}\\
&<\frac{6c_8}{\rho_0}(e^{-\frac{\rho_0}{12}2^\alpha})^{(\frac{4}{3})^n}\Big(\big((\frac{3}{2})^9d_0\big)^{(\frac{4}{3})^n}(\frac{2}{3})^{3n+9}\Big)^{-\frac{1}{3}}\\
&<\Big((\frac{6c_8}{\rho_0})^{3}e^{-\frac{\rho_0}{12}2^\alpha}\Big)^{\frac{1}{3}(\frac{4}{3})^n}\Big((\frac{3}{2})^9d_0\Big)^{-\frac{1}{3}(\frac{4}{3})^n}\big((\frac{3}{2})^{3n+9}\big)^{\frac{1}{3}}\\
&<\Big((\frac{6c_8}{\rho_0})^{3}e^{-\frac{\rho_0}{12}2^\alpha}\Big)^{\frac{1}{3}(\frac{4}{3})^n}\Big((\frac{3}{2})^9d_0\Big)^{-\frac{1}{3}(\frac{4}{3})^n}\big((\frac{3}{2})^{9}\big)^{\frac{1}{3}(\frac{4}{3})^n}\\
&=\Big((\frac{6c_8}{\rho_0})^{3}(\frac{3}{2})^{9}\frac{1}{e^{\frac{\rho_0}{12}2^\alpha}(\frac{3}{2})^9d_0}\Big)^{\frac{1}{3}(\frac{4}{3})^n}\\
&=\Big((\frac{6c_8}{\rho_0})^{3}\frac{1}{e^{\frac{\rho_0}{12}2^\alpha}d_0}\Big)^{\frac{1}{3}(\frac{4}{3})^n}.
\end{aligned}
\end{equation*}
Therefore when $(\frac{6c_8}{\rho_0})^{3}\frac{1}{e^{\frac{\rho_0}{12}2^\alpha}d_0}<1$,
i.e.,
\begin{equation}\label{esti15}
\alpha>\log_{2}^{\frac{12}{\rho_0}\ln{\frac{216c_9}{\rho_0^3d_0}}},
\end{equation}
we have
 $$c_8(\frac{\rho_0}{6})^{-1}e^{-\frac{\rho_0}{12}(\frac{2}{3})^n(2s_n-1)}d_n^{-\frac{1}{3}}<1,$$
 and
 $$\frac{c_6(\rho_n-\rho_{n+1})^{-1}d_ne^{-\frac{1}{2}(2s_n-1)(\rho_n-\rho_{n+1})}}{1-\nu}<\frac{1}{2}d_{n+1},
 $$
where $c_9=c_8^3$.

On the other hand, we will prove the inequality $$\frac{c_6(\rho_n-\rho_{n+1})^{-1}d_n^2(t_n-t_{n+1})^{-m-\sigma}(\frac{1}{t_n-t_{n+1}}+\frac{1}{\rho_n-\rho_{n+1}})}{1-\nu}<\frac{1}{2}d_{n+1}.$$
Due to
\begin{equation*}
\begin{aligned}&\frac{c_6(\rho_n-\rho_{n+1})^{-1}d_n^2(t_n-t_{n+1})^{-m-\sigma}(\frac{1}{t_n-t_{n+1}}+\frac{1}{\rho_n-\rho_{n+1}})}{1-\nu}\\
&<\frac{5}{4}c_6\frac{36}{\rho_0}(\frac{6}{t_0})^{m+\sigma}(\frac{1}{\rho_0}+\frac{1}{t_0})(\frac{3}{2})^{n(m+\sigma+1)}d_n^{\frac{2}{3}}(\frac{3}{2})^nd_n^{\frac{4}{3}}\\
&<\frac{1}{2}\Big(\frac{36c_8}{\rho_0}(\frac{6}{t_0})^{m+\sigma}(\frac{1}{\rho_0}+\frac{1}{t_0})(\frac{3}{2})^{n(m+\sigma+1)}d_n^{\frac{2}{3}}\Big)d_{n+1},
\end{aligned}
\end{equation*}
we turn to prove \begin{equation}\label{esti14}
\frac{36c_8}{\rho_0}(\frac{6}{t_0})^{m+\sigma}(\frac{1}{\rho_0}+\frac{1}{t_0})(\frac{3}{2})^{n(m+\sigma+1)}d_n^{\frac{2}{3}}<1.
\end{equation}
The estimate of the left hand in (\ref{esti14}) follows
\begin{equation*}
\begin{aligned}&\frac{36c_8}{\rho_0}(\frac{6}{t_0})^{m+\sigma}(\frac{1}{\rho_0}+\frac{1}{t_0})(\frac{3}{2})^{n(m+\sigma+1)}d_n^{\frac{2}{3}}\\
&=\frac{36c_8}{\rho_0}(\frac{6}{t_0})^{m+\sigma}(\frac{1}{\rho_0}+\frac{1}{t_0})(\frac{3}{2})^{n(m+\sigma+1)}\Big(\big((\frac{3}{2})^9d_0\big)^{(\frac{4}{3})^n}(\frac{2}{3})^{3n+9}\Big)^{\frac{2}{3}}\\
&<\frac{36c_8}{\rho_0}(\frac{6}{t_0})^{m+\sigma}(\frac{1}{\rho_0}+\frac{1}{t_0})(\frac{3}{2})^{n(m+\sigma-1)-6}\big((\frac{3}{2})^6d_0^{\frac{2}{3}}\big)^{(\frac{4}{3})^n}\\
&<\Big(\frac{36c_8}{\rho_0}(\frac{6}{t_0})^{m+\sigma}(\frac{1}{\rho_0}+\frac{1}{t_0})(\frac{3}{2})^{(m+\sigma-1)}(\frac{3}{2})^6d_0^{\frac{2}{3}}\Big)^{(\frac{4}{3})^n}\\
&=\Big(\frac{36c_8}{\rho_0}(\frac{6}{t_0})^{m+\sigma}(\frac{1}{\rho_0}+\frac{1}{t_0})(\frac{3}{2})^{(m+\sigma+5)}d_0^{\frac{2}{3}}\Big)^{(\frac{4}{3})^n}.
\end{aligned}
\end{equation*}
Thus, by the choice of
\begin{equation}\label{esti16}
d_0<[\frac{\rho_0}{36c_8}(\frac{t_0}{6})^{m+\sigma}(\frac{1}{\rho_0}+\frac{1}{t_0})^{-1}(\frac{2}{3})^{(m+\sigma+5)}]^{\frac{3}{2}},
\end{equation}
 (\ref{esti14}) is established.

From the above discussion, we conclude
$$|f_{s_{n+1}}(\xi,\eta)|_{D_{n+1}}+|g_{s_{n+1}}(\xi,\eta)|_{D_{n+1}}<d_{n+1},$$
if (\ref{esti15}) and (\ref{esti16}) hold.
Hence we get the inequality (\ref{pert-esti}), which completes the induction.

For one thing, the change of variables $T_n$ can be written as $T_n=T_0\circ\Delta T_1\circ\Delta T_2\circ\cdots\Delta T_n$, then $ T_{n+1}=T_n\circ \Delta T_{n+1}$, and the transformation $T_{n+1}$ can be expressed as
\begin{equation*}T_{n+1}:
\begin{cases}
\theta=\xi+u_{n+1}(\xi,\eta)\\
r=\eta+v_{n+1}(\xi,\eta),
\end{cases}
\end{equation*}
where
\begin{equation}\label{nonli}
u_{n+1}=u_0+u_1+...+u_n,\quad \quad
v_{n+1}=v_0+v_1+...+v_n.
\end{equation}
The convergence of the transformation sequence $\{T_{n+1}\}$ is decided by their nonlinear parts (\ref{nonli}).

By Lemma \ref{iter}, we have
\begin{equation*}
\begin{aligned}
&|u_n|_{D_n}+|v_n|_{D_n}<c_2d_n(t_n-t_{n+1})^{-m-\sigma}(\rho_n-\rho_{n+1})^{-1}\\
&=c_2(\frac{6}{t_0})^{m+\sigma}\frac{6}{\rho_0}(\frac{3}{2})^{n(m+\sigma+1)}(\frac{2}{3})^{3n+9}\big((\frac{3}{2})^{9}d_0\big)^{(\frac{4}{3})^n}.
\end{aligned}
\end{equation*}
When $d_0<(\frac{2}{3})^{m+\sigma+7}$, there is $|u_n|_{D_n}+|v_n|_{D_n}\rightarrow 0$ as $n\rightarrow\infty$. It follows that the sequence $\{u_{n+1}\}$ and $\{v_{n+1}\}$ are uniformly bounded in $D_\infty=\{(\xi,\eta):|\Im\xi|<\frac{t_0}{2}, |\eta|<\frac{\rho_0}{2}\}$. Hence, one can choose a subsequence of $\{T_{n}\}$, which converges to a transformation $T$ in $D_\infty$.

For another, since the nonlinear parts of $M_n$ satisfies $|f_n|_{D_n}+|g_n|_{D_n}<d_n\rightarrow 0$ as $n\rightarrow\infty$, it implies that the mapping $M_n={T_n}^{-1}\circ M_{0}\circ T_n$ tends to a linearized normal form (\ref{norm-form}) in $D_\infty$.

In conclusion, there is a convergent transformation $V\circ T$, where $V$ is defined in Lemma \ref{prop}, such that the mapping (\ref{main-equa}) is reduced to the formal normal form (\ref{norm-form}).

\section{Application}

 In this section, we will apply Theorem \ref{main-theo} to the equation
 \begin{equation}\label{appl-equa}
x''+\varphi(x)f(x')+\omega^2x+g(x)=p(t).
 \end{equation}

Suppose that

$(H_1)$: functions $\varphi, f, g$ and $p$ are real analytic in $x,x'$ and $t$;

$(H_2)$: $f$ and $p$ are even functions, $p(t+2\pi)=p(t)$;

$(H_3)$: $\lim_{x\rightarrow\pm\infty}\varphi(x)=:\varphi(\pm\infty)\in \mathbb{R},\quad \lim_{|x|\rightarrow+\infty}x^4\varphi^{(4)}(x)=0;$

$(H_4)$: $\lim_{x\rightarrow\pm\infty}f(x)=:f(+\infty)\in \mathbb{R},\quad \lim_{|x|\rightarrow+\infty}x^4f^{(4)}(x)=0;$

$(H_5)$: $\lim_{x\rightarrow\pm\infty}g(x)=:g(\pm\infty)\in \mathbb{R},\quad \lim_{|x|\rightarrow+\infty}x^4g^{(4)}(x)=0.$

\begin{theorem}\label{appl-theo}
Suppose that $(H_1)-(H_5)$ hold, and
$\omega$ satisfies
$$|k\omega^{-1}-l|\geq\frac{c_0}{|k|^{\sigma}},$$
where $c_0>0,\sigma>0, k\in\mathbb{Z}\setminus 0, l\in\mathbb{Z}$.
Then for every solution $x(t)$ of (\ref{appl-equa}), we have
$$\sup_{t\in\mathbb{R}}(|x(t)|+|x'(t)|)<+\infty.$$
\end{theorem}

In order to obtain the boundedness of all solutions of (\ref{appl-equa}), it is sufficient to prove that its Poincar\'{e} mapping can be written as a twist mapping with small enough perturbations. Under some transformations, if some Birkhoff constants of the Poincar\'{e} mapping are not zero, we use classical twist theorem or small twist theorem for reversible mappings to derive the boundedness. Otherwise, if all Birkhoff constants of the Poincar\'{e} mapping vanish, we apply Theorem \ref{main-theo} to achieve the goal.  In the following, we will give the proof of Theorem \ref{appl-theo}, which is similar to the proof in \cite{Liubounde} and \cite{LiM}. Thus we give the sketch of the proof.

We first rewrite (\ref{appl-equa}) as
\begin{equation}\label{appl-syst1}
\begin{cases}
x'=-\omega y\\
y'=\omega x+\omega^{-1}\varphi(x)f(\omega y)+\omega^{-1}g(x)-\omega^{-1}p(t).
\end{cases}
\end{equation}
From $(H_2)$, it follows that (\ref{appl-syst1}) is reversible with respect to the involution $G(x,y)=(x,-y)$.

By polar coordinates change
$x=r\cos\theta,\ y=r\sin\theta,$
the system (\ref{appl-syst2}) is transformed into
\begin{equation}\label{appl-syst2}
\begin{cases}
r'=\omega^{-1}\big(\varphi(r\cos\theta)f(\omega r\sin\theta)+g(r\cos\theta)\big)\sin\theta-\omega^{-1}p(t)\sin\theta\\
\theta'=\omega +\omega^{-1}r^{-1}\big(\varphi(r\cos\theta)f(\omega r\sin\theta)+g(r\cos\theta)\big)\cos\theta-\omega^{-1}r^{-1}p(t)\cos\theta.
\end{cases}
\end{equation}
Observing that
$$|\omega^{-1}r^{-1}\big(\varphi(r\cos\theta)f(\omega r\sin\theta)+g(r\cos\theta)\big)\cos\theta-\omega^{-1}r^{-1}p(t)\cos\theta|\leq Cr^{-1},$$
for some $C>0$, we may consider (\ref{appl-syst2}) assuming that $r(t)>2C\omega^{-1}$ for all $t\in\mathbb{R}$ along a solution $t\mapsto(r(t),\theta(t))$. Therefore,
$$\theta'\geq\frac{1}{2}\omega>0,\quad t\in\mathbb{R},$$
which means that $t\mapsto\theta(t)$ ia globally invertible. Denoting by $\theta\mapsto t(\theta)$ the inverse function, we have that $\theta\mapsto (r(t(\theta)),t(\theta))$ solves the system
\begin{equation}\label{appl-syst3}
\begin{cases}
\frac{dr}{d\theta}=\Phi(r,t,\theta)\\
\frac{dt}{d\theta}=\Psi(r,t,\theta),
\end{cases}
\end{equation}
where $$\Phi(r,t,\theta)=\frac{\omega^{-1}\big(\varphi(r\cos\theta)f(\omega r\sin\theta)+g(r\cos\theta)\big)\sin\theta-\omega^{-1}p(t)\sin\theta}{\omega +\omega^{-1}r^{-1}\big(\varphi(r\cos\theta)f(\omega r\sin\theta)+g(r\cos\theta)\big)\cos\theta-\omega^{-1}r^{-1}p(t)\cos\theta},$$
$$\Psi(r,t,\theta)=\frac{1}{\omega +\omega^{-1}r^{-1}\big(\varphi(r\cos\theta)f(\omega r\sin\theta)+g(r\cos\theta)\big)\cos\theta-\omega^{-1}r^{-1}p(t)\cos\theta}.$$
Now noting that the action, angle and time variables are $r, t$ and $\theta$, respectively. Since $\Psi(r,-t,-\theta)=\Psi(r,t,\theta)$ and $\Phi(r,-t,-\theta)=-\Phi(r,t,\theta)$, we see that system (\ref{appl-syst3}) is reversible under the transformation $(r,t)\mapsto(r,-t)$.

To estimate error terms, we introduce some notations.
\begin{definition}\label{defin} \item[{\textrm{(i)}}]: Assume function $f(\theta,r,t)$ is $O_n(r^{-j})$, if $f$ is smooth in $(r,t)$, continue in $\theta$, periodic of period $2\pi$ in $\theta$ and $t$, moreover
$$|r^{k+j}\frac{\partial^{k+l}f}{\partial r^k\partial t^l}|\leq C,\quad 0\leq k+l\leq n,$$
where $C$ is a positive constant.

\item[{\textrm{(ii)}}]: Suppose function $f(\theta,r,t)$ is $o_n(r^{-j})$, if $f$ is smooth in $(r,t)$, continue in $\theta$, periodic of period $2\pi$ in $\theta$ and $t$, moreover
$$\lim_{r\rightarrow\infty}r^{k+j}\frac{\partial^{k+l}f}{\partial r^k\partial t^l}=0,\quad 0\leq k+l\leq n,$$
uniformly in $(\theta,t)$.
\end{definition}
It is obvious that $$\Phi(r,t,\theta)\in O_4(1),\quad\Psi(r,t,\theta)\in O_4(1), $$
and (\ref{appl-syst3}) can be rewritten as
\begin{equation}\label{appl-syst4}
\begin{cases}
\frac{dr}{d\theta}=\omega^{-2}\big(\varphi(r\cos\theta)f(\omega r\sin\theta)+g(r\cos\theta)\big)\sin\theta-\omega^{-2}p(t)\sin\theta+O_4(r^{-1})\\
\frac{dt}{d\theta}=\omega^{-1}-\omega^{-3}r^{-1}\big(\varphi(r\cos\theta)f(\omega r\sin\theta)+g(r\cos\theta)\big)\cos\theta+\omega^{-3}r^{-1}p(t)\cos\theta+O_4(r^{-2}).
\end{cases}
\end{equation}

Since the Poincar\'{e} mapping of (\ref{appl-syst4}) is not sufficiently close to a twist map, we need to transform (\ref{appl-syst4}) further.

Let
$$\lambda=r+S_1(\theta,r),\quad t=t,$$
where
$$S_1(\theta,r)=- \omega^{-2}\int_0^\theta\big(\varphi(r\cos\phi)f(\omega r\sin\phi)+g(r\cos\phi)\big)\sin\phi d\phi.$$
Under this transformation, (\ref{appl-syst4}) is transformed into
\begin{equation}\label{appl-syst5}
\begin{cases}
\frac{d\lambda}{d\theta}=-\omega^{-2}p(t)\sin\theta+O_4(\lambda^{-1}),\\
\frac{dt}{d\theta}=\omega^{-1}-\omega^{-3}\lambda^{-1}\big(\varphi(\lambda\cos\theta)f(\omega \lambda\sin\theta)+g(\lambda\cos\theta)\big)\cos\theta+\omega^{-3}\lambda^{-1}p(t)\cos\theta+O_4(\lambda^{-2}).
\end{cases}
\end{equation}
Introduce a transformation $$\lambda=\lambda,\quad \tau=t+S_2(\theta,\lambda),$$
where $$S_2(\theta,\lambda)=\omega^{-3}\lambda^{-1}\int_0^\theta\big(\varphi(\lambda\cos\phi)f(\omega \lambda\sin\phi)\cos\phi-\lambda J_1(\lambda)\big)+\big(g(\lambda\cos\phi)\cos\phi-\lambda J_2(\lambda)\big) d\phi,$$
with $$J_1(\lambda)=\frac{1}{2\pi\lambda}\int_0^{2\pi}\varphi(\lambda\cos\phi)f(\omega \lambda\sin\phi)\cos\phi d\phi,$$
$$J_2(\lambda)=\frac{1}{2\pi\lambda}\int_0^{2\pi}g(\lambda\cos\phi)\cos\phi d\phi.$$
By this transformation, (\ref{appl-syst5}) is transformed into
\begin{equation}\label{appl-syst6}
\begin{cases}
\frac{d\lambda}{d\theta}=-\omega^{-2}p(\tau)\sin\theta+O_4(\lambda^{-1})\\
\frac{d\tau}{d\theta}=\omega^{-1}-\omega^{-3}\big(J_1(\lambda)+J_2(\lambda)\big)+\omega^{-3}\lambda^{-1}p(\tau)\cos\theta+O_4(\lambda^{-2}).
\end{cases}
\end{equation}
Furthermore, we can find a transformation $(\lambda,\tau)\rightarrow(\lambda,\varsigma)$,
where \begin{equation}\label{change2}
\varsigma=\tau+\lambda^{-1}S_3(\theta,\tau),
\end{equation}
and $S_3(\theta,\tau)$ is determined by solving equation
\begin{equation*}
\omega^{-3}p(\tau)\cos\theta+\frac{\partial S_3}{\partial\theta}+\omega^{-1}\frac{\partial S_3}{\partial\tau}=0.
\end{equation*}
By this transformation, we eliminate $\omega^{-3}\lambda^{-1}p(\tau)\cos\theta$ item in the second equation of (\ref{appl-syst6}).

Since $$\lim_{\lambda\rightarrow+\infty}\lambda^{k+1}J_1^{(k)}(\lambda)=(-1)^kk!\frac{1}{\pi}\big(\varphi(+\infty)-\varphi(-\infty)\big)f(+\infty),$$ and
$$\lim_{\lambda\rightarrow+\infty}\lambda^{k+1}J_2^{(k)}(\lambda)=(-1)^kk!\frac{1}{\pi}(g(+\infty)-g(-\infty)), \quad0\leq k\leq 4,$$
 (\ref{appl-syst6}) can be rewritten as
\begin{equation}\label{appl-syst7}
\begin{cases}
\frac{d\lambda}{d\theta}=-\omega^{-2}p(\varsigma)\sin\theta+O_4(\lambda^{-1})\\
\frac{d\varsigma}{d\theta}=\omega^{-1}-\frac{\omega^{-3}}{\pi}\lambda^{-1}\Big(\big(\varphi(+\infty)-\varphi(-\infty)\big)f(+\infty)+\big(g(+\infty)-g(-\infty)\big)\Big)+o_4(\lambda^{-1}).
\end{cases}
\end{equation}
We also recalled (\ref{appl-syst7}) is reversible with respect to $(\lambda,\varsigma)\mapsto(\lambda,-\varsigma)$.

Denote $\lambda=\rho^{-1},$ then (\ref{appl-syst7}) can be rewritten as

\begin{equation*}
\begin{cases}
\frac{d\rho}{d\theta}=\omega^{-2}\rho^2p(\varsigma)\sin\theta+o_4(\rho^2),\\
\frac{d\varsigma}{d\theta}=\omega^{-1}-\frac{\omega^{-3}}{\pi}\rho\Big(\big(\varphi(+\infty)-\varphi(-\infty)\big)f(+\infty)+\big(g(+\infty)-g(-\infty)\big)\Big)+o_4(\rho).
\end{cases}
\end{equation*}
Therefore, we derive the corresponding Poincar\'{e} mapping with the form
 \begin{equation}\label{appl-poin}
\begin{cases}
\rho_1=\rho_0+\rho_0^2 l(\varsigma_0)+o_4(\rho_0^2),\\
\varsigma_1=\varsigma_0+\gamma_0+\gamma_1\rho_0+o_4(\rho_0),
\end{cases}
\end{equation}
where $\gamma_0=2\pi\omega^{-1}$, $l(\varsigma_0)=\omega^{-2}\int_0^{2\pi} p(\varsigma_0+\omega^{-1}\theta)\sin\theta d\theta$
and
 \begin{equation}\label{twistcondition}
 \gamma_1=-2\omega^{-3}\big((\varphi(+\infty)-\varphi(-\infty))f(+\infty)+(g(+\infty)-g(-\infty))\big).
 \end{equation}
From the normal form theory, we see that one of the following two cases occurs.

Case 1. If $(\varphi(+\infty)-\varphi(-\infty))f(+\infty)+(g(+\infty)-g(-\infty))\neq0$, that is, $\gamma_1\neq0$, then by twist theorem for reversible mappings (see \cite{Liu-1}), there are many invariant curves for $\rho\ll 1$. If $\gamma_1=0$, we need to continue looking for changes, such that (\ref{appl-poin}) is of the following form
 \begin{equation*}
\begin{cases}
\rho_1=\rho_0+c_1\rho_0^k+o(\rho_0^{k+1}),\\
\varsigma_1=\varsigma_0+\gamma_0+\gamma_l\rho_0^l+o(\rho_0^{l+1}),
\end{cases}
\end{equation*}
where $c_1,\gamma_1$ are constants and $l>k$. Similarly, as long as the coefficient of the twist term $\gamma_l$ is not zero, by twist theorem in \cite{Liu-1}, we obtain many invariant curves if $\rho\ll 1$, i.e., $r\gg 1$, which
implies the existence of quasi-periodic solutions of (\ref{appl-equa}).

Case 2. There is a change of variables such that (\ref{appl-poin}) can be transformed into a linear normal form
 \begin{equation*}
\begin{cases}
\rho_1=\rho_0\\
\varsigma_1=\varsigma_0+\gamma_0.
\end{cases}
\end{equation*}
By Theorem \ref{main-theo}, there are sequence of invariant curves tending to $\rho_0=0$.

In conclusion, the mapping (\ref{appl-poin}) has many invariant curves tending to $\rho_0=0$, which means the invariant curves of the Poincar\'{e} map of (\ref{appl-equa}) tend to infinity. Thus for equation (\ref{appl-equa}), the existence of quasi-periodic solutions is got. Moreover for the initial value lying between two invariant curves, the solution is globally bounded. As the invariant curves tend to infinity, all solutions of (\ref{appl-equa}) are bounded, therefore we finish the proof of Theorem \ref{appl-theo}.

\section*{Acknowledgment}

Y. Z was partially supported by the NSFC (grant no. 11971059) and the Fundamental Research Funds for the Central Universities (grant no. 202261096).
D. P was supported in part by the NSFC (grant no. 11971059).  Y. N was supported by the NSFC (grant no. 12201587) , the China
Postdoctoral Science Foundation (grant no. 2020M682236) and the Shandong Provincial Natural Science Foundation, China (grant no. A010704).



\end{spacing}
\end{document}